\documentclass[oneside,english,secthm,seceqn]{elsart1p}
\usepackage{mathptmx}
\usepackage[T1]{fontenc}
\usepackage[latin9]{inputenc}
\usepackage{amsmath}
\usepackage{amssymb}

\makeatletter
\newtheorem{theorema}{Theorem A}

\usepackage[verbose]{cite}
\usepackage{babel}
\usepackage{mathptmx}

\makeatother

\journal{Journal of Difference Equations and Applications}

\makeatother

\makeatother

\makeatother

\makeatother

\begin{document}

\begin{frontmatter}
\title{On a generalization of the global attractivity for a periodically
forced Pielou's equation}
\author[waseda]{Keigo Ishihara}
\ead{keigo.i.123@gmail.com}
\author[waseda]{Yukihiko Nakata\corauthref{cor1}}
\ead{yunayuna.na@gmail.com}
\corauth[cor1]{Corresponding author}
\address[waseda]{Department of Pure and Applied Mathematics, Waseda University, 3-4-1
Ohkubo, Shinjuku-ku, Tokyo 169-8555, Japan}
\begin{abstract}
In this paper, we study the global attractivity for a class of periodic
difference equation with delay which has a generalized form of Pielou's
difference equation. The global dynamics of the equation is characterized
by using a relation between the upper and lower limit of the solution.
There are two possible global dynamics: zero solution is globally
attractive or there exists a periodic solution which is globally attractive.
Recent results in {[}E. Camouzis, G. Ladas, Periodically forced Pielou's
equation, J. Math. Anal. Appl. 333 (1) (2007) 117-127{]} is generalized.
Two examples are given to illustrate our results.\end{abstract}
\end{frontmatter}

\section{Introduction }

Several authors have studied difference equations for mathematical
models in population biology (see \cite{MR0252051,MR1162273,MR1247956,MR2100727,MR2151684,MR2314577,MR2397325,MR2323480,MR2413754,MR2447183,Liz200626}
and the references therein). Pielou's equation was proposed by Pielou
in \cite{MR0252051} as a discrete analogue of the logistic equation
with delay.

Camouzis and Ladas \cite{MR2323480} studied the following Pielou's
equation with a periodic coefficient, \begin{equation}
x_{n+1}=\frac{\beta_{n}x_{n}}{1+x_{n-1}},n=0,1,2,\dots,\label{periodicpielou}\end{equation}
where $\beta_{n},n=0,1,2,\dots$ is a periodic sequence with an arbitrary
positive integer $k$. Initial conditions are given by $x_{0}>0,x_{-1}\geq0$.
They proved that every solution converges to $0$, if $\Pi_{n=1}^{k}\beta_{n}\leq1$,
while there exists a $k$-periodic solution which is globally attractive,
if $\Pi_{n=1}^{k}\beta_{n}>1$ by using an interesting relation between
the upper and lower limit of the solution. Recently, Nyerges \cite{MR2413754}
studied the global dynamics of a general autonomous difference equation
by extending their idea. Let us introduce the main result in Camouzis
and Ladas \cite{MR2323480}.

\begin{theorema} (see Camouzis and Ladas \cite[Theorems 3.2, 3.3 and 3.4]{MR2323480})
If \[
\Pi_{n=1}^{k}\beta_{n}>1,\]
then there exists a $k$-periodic solution $x_{n}^{*}$ such that
$x_{n}^{*}=x_{n+k}^{*}$ which is globally attractive, that is, for
any solution of (\ref{periodicpielou}), it holds that \[
\lim_{n\to+\infty}(x_{n}-x_{n}^{*})=0.\]
\end{theorema}

In this paper, we further generalize Theorem A. The present paper
is focused on nonautonomous difference equations, different from Nyerges
\cite{MR2413754}. We shall study the following difference equation,
\begin{equation}
x_{n+1}=x_{n}f_{n}(x_{n-1}),n=0,1,2,\dots,\label{eq:Mother_generalPielou}\end{equation}
where $f_{n}(x),n=0,1,2,\dots,$ is continuous, bounded and positive
function on $x\in[0,+\infty)$ and $k$-periodic on $n$, that is,\[
f_{n}(x)=f_{n+k}(x)\text{ where }k\text{ is an arbitrary positive integer}.\]
It is assumed that the initial conditions are given by $x_{0}=x^{0}>0,x_{-1}=x^{-1}\geq0$.
For the function $f_{n}(x),n=1,2,\dots,k$, we impose the following
monotonicity property\begin{equation}
f_{n}(x)\text{ is strictly decreasing on }x\in[0,+\infty)\text{ for }n=1,2,\dots,k,\label{eq:ass00}\end{equation}
and\begin{equation}
xf_{n}(x)\text{ is strictly increasing on }x\in[0,+\infty)\text{ for }n=1,2,\dots,k.\label{eq:assumption}\end{equation}
Under the assumption (\ref{eq:ass00}), $\Pi_{n=1}^{k}f_{n}(x)$ is
also strictly decreasing on $x\in[0,+\infty)$ and hence, there exists
a some constant $c\geq0$ such that \[
c=\lim_{x\to+\infty}\Pi_{n=1}^{k}f_{n}(x).\]
The following theorem is a generalized version of Theorem A.
\begin{thm}
\label{thm:mainthm}Assume that (\ref{eq:ass00}) and (\ref{eq:assumption}).
If \begin{equation}
\Pi_{n=1}^{k}f_{n}(0)>1,\text{ and }\lim_{x\to+\infty}\Pi_{n=1}^{k}f_{n}(x)<1,\label{eq:theorem11}\end{equation}
 then there exists a $k$-periodic solution $x_{n}^{*}$ such that
$x_{n}^{*}=x_{n+k}^{*}$ which is globally attractive, that is, for
any solution of (\ref{eq:Mother_generalPielou}), it holds that \[
\lim_{n\to+\infty}(x_{n}-x_{n}^{*})=0.\]

\end{thm}
One can see that the assumptions (\ref{eq:ass00}) and (\ref{eq:assumption})
are nice properties to obtain the global character of the periodic
difference equation which has a form (\ref{eq:Mother_generalPielou}).
Motivated by Camouzis and Ladas \cite{MR2323480}, we also find a
useful relation between the upper and lower limit of the solution.
Then we can obtain the existence of $k$-periodic solution which is
globally attractive if (\ref{eq:theorem11}) holds. 

The paper is organized as follows. At first, we consider the case
where every solution approaches to zero solution in Section 2. In
Section 3, we show that every solution is bounded above and below
by a positive constant, respectively. This makes possible to consider
a set of the upper and lower limit of the solutions which are positive
constants (see (\ref{eq:defsh})). In Sections 4 and 5, we consider
the existence of a $k$-periodic solution which is globally attractive.
We divide the discussion in two cases, $k$ is an even integer in
Section 4 and $k$ is an odd integer in Section 5. It is important
to establish the relation (see Lemmas \ref{lem:Lemma_even} and \ref{lem:Lemma_odd})
between the upper and lower limit of the solution in these sections.
In Section 6, we apply Theorems \ref{thm:mainthm} to two nonautonomous
difference equations. The global attractivity for a delayed Beverton-Holt
equation with a periodic coefficient is established.

\section{Global attractivity of zero solution}

First of all, we consider the case where every solution approaches
to zero. Let us introduce the following result which generalizes Theorem
3.1 in Camouzis and Ladas \cite{MR2323480}.
\begin{thm}
\label{thm:main0}Assume that (\ref{eq:ass00}). If\begin{equation}
\Pi_{n=1}^{k}f_{n}(0)\leq1,\label{eq:theorem12}\end{equation}
then, for any solution of (\ref{eq:Mother_generalPielou}), it holds
that\[
\lim_{n\to+\infty}x_{n}=0.\]
\end{thm}
\begin{pf}
We have \[
\Pi_{n=1}^{k}f_{n}(x_{n})\begin{cases}
<\Pi_{n=1}^{k}f_{n}(0) & \text{ if there exists \ensuremath{n\in\left\{ 1,2,\dots,k\right\} } such that \ensuremath{x_{n}>0,}}\\
=\Pi_{n=1}^{k}f_{n}(0) & \text{ for }x_{n}=0,n=1,2,\dots,k,\end{cases}\]
from (\ref{eq:ass00}). Since $x_{n}>0$ for $n=1,2,\dots,$ for given
the initial conditions $x_{0}=x^{0}>0,x_{-1}=x^{-1}\geq0$, it holds
\[
x_{n}=x_{n-k}\Pi_{j=1}^{k}f_{n-j}(x_{n-j-1})<x_{n-k},\text{ for }n=k+1,k+2,\dots.\]
Thus, we obtain the conclusion and the proof is complete.\qed 
\end{pf}

\section{Permanence}

In this section, we show that every solution is bounded above and
below by a positive constant, respectively, if (\ref{eq:theorem11})
holds. Based on the following result, we investigate the existence
of the periodic solution which is globally attractive for any solution
in Sections 4 and 5.
\begin{thm}
\label{thm:Permanence}Assume that (\ref{eq:ass00}). If (\ref{eq:theorem11}),
then, for any solution of (\ref{eq:Mother_generalPielou}), it holds
that\[
0<\tilde{x}\Pi_{n=1}^{k}f_{n}(\overline{x})\leq\liminf_{n\to+\infty}x_{n}\leq\limsup_{n\to+\infty}x_{n}\leq\tilde{x}\Pi_{n=1}^{k}f_{n}(0)<+\infty,\]
where $\overline{x}=\tilde{x}\Pi_{n=1}^{k}f_{n}(0)$ and $\tilde{x}$
is a unique positive solution of $\Pi_{n=1}^{k}f_{n}(x)=1$.\end{thm}
\begin{pf}
We see that there exists $\tilde{x}<+\infty$ such that $\Pi_{n=1}^{k}f_{n}(\tilde{x})=1$
by (\ref{eq:theorem11}).

Suppose that $\limsup_{n\to+\infty}x_{n}=+\infty$. Then there exists
a subsequence $\{\overline{n}_{m}\}_{m=1}^{+\infty}$ such that\begin{equation}
x_{\overline{n}_{m}}=\max_{0\leq n\leq\overline{n}_{m}}x_{n}\text{ and }\lim_{m\to+\infty}x_{\overline{n}_{m}}=\limsup_{n\to+\infty}x_{n}=+\infty.\label{eq:suppose1}\end{equation}
From (\ref{eq:Mother_generalPielou}), it holds that\begin{align}
x_{\overline{n}_{m}} & =x_{\overline{n}_{m}-1}f_{\overline{n}_{m}-1}(x_{\overline{n}_{m}-2})\nonumber \\
 & =x_{\overline{n}_{m}-2}f_{\overline{n}_{m}-2}(x_{\overline{n}_{m}-3})f_{\overline{n}_{m}-1}(x_{\overline{n}_{m}-2})\nonumber \\
 & =\dots\nonumber \\
 & =x_{\overline{n}_{m}-k}\Pi_{j=1}^{k}f_{\overline{n}_{m}-j}(x_{\overline{n}_{m}-j-1}).\label{eq:ite1}\end{align}
From (\ref{eq:suppose1}) and (\ref{eq:ite1}), we see\[
\Pi_{j=1}^{k}f_{\overline{n}_{m}-j}(x_{\overline{n}_{m}-j-1})\geq1,\]
and, by (\ref{eq:ass00}), it follows \[
\Pi_{j=1}^{k}f_{\overline{n}_{m}-j}(x_{\overline{n}_{m}-\underline{j}-1})\geq\Pi_{j=1}^{k}f_{\overline{n}_{m}-j}(x_{\overline{n}_{m}-j-1})\geq1,\]
where \[
x_{\overline{n}_{m}-\underline{j}-1}=\min_{1\leq j\le k}x_{\overline{n}_{m}-j-1}.\]
Then, we see $x_{\overline{n}_{m}-\underline{j}-1}\leq\tilde{x}$
and hence, from (\ref{eq:ite1}), we obtain \begin{align*}
x_{\overline{n}_{m}} & =x_{\overline{n}_{m}-\underline{j}-1}\Pi_{j=1}^{\underline{j}+1}f_{\overline{n}_{m}-j}(x_{\overline{n}_{m}-j-1})\\
 & \leq\tilde{x}\Pi_{n=1}^{k}f_{n}(0)<+\infty,\end{align*}
which leads a contradiction to our assumption. Thus \[
\limsup_{n\to+\infty}x_{n}<+\infty.\]

Moreover, similar to the above discussion, we obtain that\[
\limsup_{n\to+\infty}x_{n}\leq\tilde{x}\Pi_{n=1}^{k}f_{n}(0)<+\infty.\]

Next, suppose that $\liminf_{n\to+\infty}x_{n}=0$. Then there exists
a subsequence $\{\underline{n}_{m}\}_{m=1}^{+\infty}$ such that \begin{equation}
x_{\underline{n}_{m}}=\min_{0\leq n\leq\underline{n}_{m}}x_{n}\text{ and }\lim_{m\to+\infty}x_{\underline{n}_{m}}=\liminf_{n\to+\infty}x_{n}=0.\label{eq:suppose2}\end{equation}
From (\ref{eq:Mother_generalPielou}), it holds that\begin{align}
x_{\underline{n}_{m}} & =x_{\underline{n}_{m}-1}f_{\underline{n}_{m}-1}(x_{\underline{n}_{m}-2})\nonumber \\
 & =x_{\underline{n}_{m}-2}f_{\underline{n}_{m}-2}(x_{\underline{n}_{m}-3})f_{\underline{n}_{m}-1}(x_{\underline{n}_{m}-2})\nonumber \\
 & =\dots\nonumber \\
 & =x_{\underline{n}_{m}-k}\Pi_{j=1}^{k}f_{\underline{n}_{m}-j}(x_{\underline{n}_{m}-j-1}).\label{eq:ite2}\end{align}
 From (\ref{eq:suppose2}) and (\ref{eq:ite2}), we see \[
\Pi_{j=1}^{k}f_{\underline{n}_{m}-j}(x_{\underline{n}_{m}-j-1})\leq1,\]
and, by (\ref{eq:ass00}), it follows \[
\Pi_{j=1}^{k}f_{\underline{n}_{m}-j}(x_{\underline{n}_{m}-\overline{j}-1})\leq\Pi_{j=1}^{k}f_{\underline{n}_{m}-j}(x_{\underline{n}_{m}-j-1})\leq1,\]
where \[
x_{\underline{n}_{m}-\overline{j}-1}=\max_{1\leq j\le k}x_{\underline{n}_{m}-j-1}.\]
Then, we see $x_{\underline{n}_{m}-\overline{j}-1}\geq\tilde{x}$.
Hence, from (\ref{eq:ite2}), we obtain \begin{align*}
x_{\underline{n}_{m}} & \geq\tilde{x}\Pi_{n=1}^{k}f_{n}(\overline{x})>0,\end{align*}
where $\overline{x}=\tilde{x}\Pi_{n=1}^{k}f_{n}(0)$. This leads a
contradiction to our assumption. Thus, we obtain \[
\liminf_{n\to+\infty}x_{n}>0.\]

Moreover, similar to the above discussion, we obtain \[
\liminf_{n\to+\infty}x_{n}\geq\tilde{x}\Pi_{n=1}^{k}f_{n}(\overline{x})>0.\]
Hence, the proof is complete.\qed \end{pf}
\begin{rem}
The assumption (\ref{eq:assumption}) is not needed for the permanence
of the solution.
\end{rem}
Hereafter, we assume (\ref{eq:ass00}) and (\ref{eq:assumption})
hold. We are interested in the existence of $k$-periodic solution
which is globally attractive. Let us introduce some notations which
are used throughout the paper. At first, we set \begin{equation}
\begin{cases}
S_{h} & =\limsup_{n\to+\infty}x_{kn+h},\\
I_{h} & =\liminf_{n\to+\infty}x_{kn+h},\end{cases}\text{ for }h=\dots,-1,0,1,\dots,\label{eq:defsh}\end{equation}
and will show \[
S_{h}=I_{h},\text{ for }h=1,2,\dots,k,\]
holds in order to establish the existence of the periodic solution
which is globally attractive. Obviously, by Theorem \ref{thm:Permanence},
we see $(S_{h},I_{h})\in(0,\infty)\times(0,\infty)$ if (\ref{eq:theorem11})
holds. Further, from (\ref{eq:defsh}), it follows \begin{equation}
\begin{cases}
S_{h}=\limsup_{n\to+\infty}x_{kn+h} & =\limsup_{n\to+\infty}x_{k\left(n+j\right)+h}=S_{kj+h},\\
I_{h}=\liminf_{n\to+\infty}x_{kn+h} & =\liminf_{n\to+\infty}x_{k\left(n+j\right)+h}=I_{kj+h},\end{cases}\label{periodicityofsi}\end{equation}
where $j$ is a some integer and the relation (\ref{periodicityofsi})
will be used if necessary. We, further, introduce two sets of subsequences
$\{\overline{n}_{m}^{h}\}_{m=0}^{+\infty}$ and $\{\underline{n}_{m}^{h}\}_{m=0}^{+\infty}$
for $h=1,2,\dots,k$ such that\[
\begin{cases}
\lim_{m\to+\infty}x_{k\overline{n}_{m}^{h}+h} & =\limsup_{n\to+\infty}x_{kn+h}=S_{h},\\
\lim_{m\to+\infty}x_{k\underline{n}_{m}^{h}+h} & =\liminf_{n\to+\infty}x_{kn+h}=I_{h}.\end{cases}\text{ for }h=1,2,\dots,k.\]
Finally, for simplicity of the proof, we also define $f_{h}(x)$ as
\[
f_{h}(x)=f_{k+h}(x)\text{ for }h=-1,-2,\dots,-k.\]

\section{Global attractivity for the case where $k$ is an even integer}

In this section, we show that Theorem \ref{thm:mainthm} holds when
$k$ is an even integer. For the reader, we first consider the case
$k=2$ in Section 4.1. (\ref{eq:concl_Lem_2}) in Lemma \ref{lem:Lemma_2}
has an important role in this subsection. Then, we give Theorem \ref{thm:Theorem_2}
which shows that there exists a $2$-periodic solution which is globally
attractive. In Section 4.2, we generalize these results to the case
where $k$ is an arbitrary even integer.

\subsection{Case: $k=2$}

We introduce the following lemma which plays a crucial role in the
proof of Theorem \ref{thm:Theorem_2}. 
\begin{lem}
\label{lem:Lemma_2} Let $k=2$. Assume that (\ref{eq:ass00}) and
(\ref{eq:assumption}). If (\ref{eq:theorem11}), then it holds that
\begin{equation}
f_{h}(I_{h-1})f_{h-1}(S_{h})=1,\label{eq:concl_Lem_2}\end{equation}
for $h=1,2$. \end{lem}
\begin{pf}
At first, we see that \begin{align}
x_{k\overline{n}_{m}^{h}+h} & =x_{k\overline{n}_{m}^{h}+h-1}f_{k\overline{n}_{m}^{h}+h-1}(x_{k\overline{n}_{m}^{h}+h-2})\nonumber \\
 & =x_{k\overline{n}_{m}^{h}+h-2}f_{k\overline{n}_{m}^{h}+h-2}(x_{k\overline{n}_{m}^{h}+h-3})f_{k\overline{n}_{m}^{h}+h-1}(x_{k\overline{n}_{m}^{h}+h-2})\nonumber \\
 & =x_{k\overline{n}_{m}^{h}+h-2}f_{h-2}(x_{k\overline{n}_{m}^{h}+h-3})f_{h-1}(x_{k\overline{n}_{m}^{h}+h-2}),\label{eq:lemma2_ite}\end{align}
 for $h=1,2$, from (\ref{eq:Mother_generalPielou}) and by considering
the limiting equation of (\ref{eq:lemma2_ite}) (letting $m\to+\infty$)
and using (\ref{eq:ass00}) and (\ref{eq:assumption}), it follows
\begin{align}
S_{h} & \leq S_{h-2}f_{h-2}(I_{h-3})f_{h-1}(S_{h-2})=S_{h}f_{h}(I_{h-1})f_{h-1}(S_{h}).\label{eq:lemma2Sh-1}\end{align}
Similarly, (by considering the limiting equation of (\ref{eq:lemma2_ite})
for the subsequences $\{\underline{n}_{m}^{h}\}_{m=0}^{+\infty},h=1,2$
and using (\ref{eq:ass00}) and (\ref{eq:assumption})), we can obtain
\begin{equation}
I_{h}\geq I_{h-2}f_{h-2}(S_{h-3})f_{h-1}(I_{h-2})=I_{h}f_{h}(S_{h-1})f_{h-1}(I_{h}).\label{eq:lemma2Ih}\end{equation}

Consequently, by (\ref{eq:lemma2Sh-1}) and (\ref{eq:lemma2Ih}),
the following holds \[
\begin{cases}
1 & \leq f_{h}(I_{h-1})f_{h-1}(S_{h}),\\
1 & \geq f_{h}(S_{h-1})f_{h-1}(I_{h}),\end{cases}\text{ for }h=1,2,\]
from which, we obtain (\ref{eq:concl_Lem_2}) and the proof is complete.
\qed
\end{pf}
Then, we obtain the following result.
\begin{thm}
\label{thm:Theorem_2}Let $k=2$. Assume that (\ref{eq:ass00}) and
(\ref{eq:assumption}). If (\ref{eq:theorem11}), then there exists
a $2$-periodic solution $x_{n}^{*}$ such that $x_{n}^{*}=x_{n+2}^{*}$
which is globally attractive, that is, for any solution of (\ref{eq:Mother_generalPielou}),
it holds that\[
\lim_{n\to+\infty}(x_{n}-x_{n}^{*})=0.\]
\end{thm}
\begin{pf}
In order to obtain the conclusion, we will show \begin{equation}
S_{h}=I_{h},h=1,2.\label{eq:concl2}\end{equation}
 Let\[
U_{h,h-j}=\lim_{m\to+\infty}x_{k\overline{n}_{m}^{h}+h-j}\text{ for }j=1,2,\dots,5,\]
and we claim

\begin{equation}
U_{h,h-j}=\begin{cases}
S_{h-j} & \text{ for }j=2,4,\\
I_{h-j} & \text{ for }j=1,3,5,\end{cases}\label{eq:claim01}\end{equation}
for $h=1,2$. 

At first, we see that it holds \begin{align}
x_{k\overline{n}_{m}^{h}+h-j} & =x_{k\overline{n}_{m}^{h}+h-j-2}f_{k\overline{n}_{m}^{h}+h-j-2}(x_{k\overline{n}_{m}^{h}+h-j-3})f_{k\overline{n}_{m}^{h}+h-j-1}(x_{k\overline{n}_{m}^{h}+h-j-2})\nonumber \\
 & =x_{k\overline{n}_{m}^{h}+h-j-2}f_{h-j-2}(x_{k\overline{n}_{m}^{h}+h-j-3})f_{h-j-1}(x_{k\overline{n}_{m}^{h}+h-j-2}),\label{eq:ite2h-j}\end{align}
for $h=1,2$ and $j=0,1,2,\dots$. from (\ref{eq:Mother_generalPielou}).

Firstly, we show \begin{equation}
U_{h,h-j}=\begin{cases}
S_{h-j}=S_{h} & \text{ for }j=2,\\
I_{h-j}=I_{h-1} & \text{ for }j=3,\end{cases}\label{eq:U1}\end{equation}
for $h=1,2$. Suppose that there exists some $\overline{h}\in\left\{ 1,2\right\} $
such that $U_{\overline{h},\overline{h}-2}<S_{\overline{h}}$ or $U_{\overline{h},\overline{h}-3}>I_{\overline{h}-1}$.
By considering the limiting equation of (\ref{eq:ite2h-j}) with $j=0$,
it follows\[
S_{h}=U_{h,h-2}f_{h}(U_{h,h-3})f_{h-1}(U_{h,h-2}).\]
Then, by (\ref{eq:ass00}) and (\ref{eq:assumption}), it follows
$S_{\overline{h}}<S_{\overline{h}}f_{\overline{h}}(I_{\overline{h}-1})f_{\overline{h}-1}(S_{\overline{h}})$,
which implies \[
1<f_{h}(I_{h-1})f_{h-1}(S_{h}),\text{ for }h=\overline{h}.\]
This leads a contradiction to (\ref{eq:concl_Lem_2}) in Lemma \ref{lem:Lemma_2}.
Thus, (\ref{eq:U1}) holds.

Next, we show \begin{equation}
U_{h,h-j}=\begin{cases}
S_{h-j}=S_{h} & \text{ for }j=4,\\
I_{h-j}=I_{h-1} & \text{ for }j=5,\end{cases}\label{eq:U2}\end{equation}
for $h=1,2$. Suppose that there exists some $\overline{h}\in\left\{ 1,2\right\} $
such that $U_{\overline{h},\overline{h}-4}<S_{\overline{h}}$ or $U_{\overline{h},\overline{h}-5}>I_{\overline{h}-1}$.
By considering the limiting equation of (\ref{eq:ite2h-j}) with $j=2$
and substituting (\ref{eq:U1}), it follows\[
U_{h,h-2}=S_{h}=U_{h,h-4}f_{h}(U_{h,h-5})f_{h-1}(U_{h,h-4}),\]
for $h=1,2$. Then, by (\ref{eq:ass00}) and (\ref{eq:assumption}),
it follows $S_{\overline{h}}<S_{\overline{h}}f_{\overline{h}}(I_{\overline{h}-1})f_{\overline{h}-1}(S_{\overline{h}})$
which implies\[
1<f_{h}(I_{h-1})f_{h-1}(S_{h}),\text{ for }h=\overline{h}.\]
This gives a contradiction to (\ref{eq:concl_Lem_2}) in Lemma \ref{lem:Lemma_2}.
Thus, (\ref{eq:U2}) holds.

By considering the limiting equation of (\ref{eq:ite2h-j}) with $j=1$
and using (\ref{eq:U1})-(\ref{eq:U2}), we see\[
U_{h,h-1}=U_{h,h-3}f_{h-1}(U_{h,h-4})f_{h}(U_{h,h-3})=I_{h-1}f_{h-1}(S_{h})f_{h}(I_{h-1}).\]
Hence, it follows \begin{equation}
U_{h,h-1}=I_{h-1},\label{eq:U4}\end{equation}
for $h=1,2$, by (\ref{eq:concl_Lem_2}) in Lemma \ref{lem:Lemma_2}.
Consequently, (\ref{eq:claim01}) holds from (\ref{eq:U1}), (\ref{eq:U2})
and (\ref{eq:U4}).

From (\ref{eq:Mother_generalPielou}), it holds \[
x_{k\overline{n}_{m}^{h}+h}=x_{k\overline{n}_{m}^{h}+h-1}f_{k\overline{n}_{m}^{h}+h-1}(x_{k\overline{n}_{m}^{h}+h-2})=x_{k\overline{n}_{m}^{h}+h-1}f_{h-1}(x_{k\overline{n}_{m}^{h}+h-2}),\]
and by considering the limiting equation and using (\ref{eq:claim01}),
we obtain\[
S_{h}=U_{h,h-1}f_{h-1}(U_{h,h-2})=I_{h-1}f_{h-1}(S_{h}).\]
By (\ref{eq:concl_Lem_2}) in Lemma \ref{lem:Lemma_2}, it holds \begin{equation}
S_{h}f_{h}(I_{h-1})=I_{h-1}f_{h-1}(S_{h})f_{h}(I_{h-1})=I_{h-1},\label{eq2s}\end{equation}
for $h=1,2$, On the other hand, similar to the above discussion,
it also holds \begin{equation}
I_{h}f_{h}(S_{h-1})=S_{h-1}f_{h-1}(I_{h})f_{h}(S_{h-1})=S_{h-1},\label{eq2i}\end{equation}
for $h=1,2$. Consequently, it holds \[
I_{h-1}=S_{h}f_{h}(I_{h-1})\leq I_{h}f_{h}(S_{h-1})=S_{h-1},\]
by (\ref{eq2s}) and (\ref{eq2i}), and hence, (\ref{eq:concl2})
holds. Then, from (\ref{eq:defsh}), we see \[
\liminf_{n\to+\infty}x_{kn+h}=\limsup_{n\to+\infty}x_{kn+h}=\lim_{n\to+\infty}x_{kn+h},\]
for $h=1,2$ and there exist two positive constants $x_{1}^{*}=S_{1}=I_{1}$
and $x_{2}^{*}=S_{2}=I_{2}$ such that \[
x_{2}^{*}=\lim_{n\to+\infty}x_{2n},\]
and \[
x_{1}^{*}=\lim_{n\to+\infty}x_{2n+1}.\]
The proof is complete.\qed 
\end{pf}

\subsection{Case: $k$ is an even integer}

In this subsection, we generalize results in Section 4.1 to the case
where $k$ is an arbitrary even integer. 
\begin{lem}
\label{lem:Lemma_even}Let $k$ be an even integer. Assume that (\ref{eq:ass00})
and (\ref{eq:assumption}). If (\ref{eq:theorem11}), then it holds
that\begin{equation}
\left(f_{h-k}(I_{h-k-1})f_{h+1-k}(S_{h-k})\right)\dots\left(f_{h-4}(I_{h-5})f_{h-3}(S_{h-4})\right)\left(f_{h-2}(I_{h-3})f_{h-1}(S_{h-2})\right)=1,\label{eq:concl_lem_even}\end{equation}
 for $h=1,2,\dots,k.$ \end{lem}
\begin{pf}
At first, from (\ref{eq:Mother_generalPielou}), we see that it holds\begin{align}
x_{k\overline{n}_{m}^{h}+h} & =x_{k\overline{n}_{m}^{h}+h-2}f_{k\overline{n}_{m}^{h}+h-2}(x_{k\overline{n}_{m}^{h}+h-3})f_{k\overline{n}_{m}^{h}+h-1}(x_{k\overline{n}_{m}^{h}+h-2})\nonumber \\
 & =x_{k\overline{n}_{m}^{h}+h-2}f_{h-2}(x_{k\overline{n}_{m}^{h}+h-3})f_{h-1}(x_{k\overline{n}_{m}^{h}+h-2}),\label{eq:lemma3_ite}\end{align}
 for $h=1,2,\dots,k$, and by considering the limiting equation of
(\ref{eq:lemma3_ite}) and using (\ref{eq:ass00}) and (\ref{eq:assumption}),
it follows\begin{equation}
S_{h}\leq S_{h-2}f_{h-2}(I_{h-3})f_{h-1}(S_{h-2}),\label{eq:h-j}\end{equation}
for $h=1,2,\dots,k$. From (\ref{eq:h-j}) we obtain\begin{align*}
S_{h} & \leq S_{h-2}\left(f_{h-2}(I_{h-3})f_{h-1}(S_{h-2})\right)\\
 & \leq S_{h-4}\left(f_{h-4}(I_{h-5})f_{h-3}(S_{h-4})\right)\left(f_{h-2}(I_{h-3})f_{h-1}(S_{h-2})\right)\\
 & \leq\dots\\
 & \leq S_{h-k}\overline{G}_{h},\end{align*}
where \begin{align}
\overline{G}_{h} & =\left(f_{h-k}(I_{h-k-1})f_{h+1-k}(S_{h-k})\right)\dots\left(f_{h-4}(I_{h-5})f_{h-3}(S_{h-4})\right)\left(f_{h-2}(I_{h-3})f_{h-1}(S_{h-2})\right).\label{eq:overg}\end{align}
Then \[
1\leq\overline{G}_{h},\]
holds for $h=1,2,\dots,k$, because we have $S_{h}=S_{h-k}$ from
(\ref{periodicityofsi}).

Similarly (by considering the limiting equation of (\ref{eq:lemma3_ite})
for the subsequences $\{\underline{n}_{m}^{h}\}_{m=0}^{+\infty},h=1,2,\dots,k$
and using (\ref{eq:ass00}) and (\ref{eq:assumption})), it follows\begin{equation}
I_{h}\geq I_{h-2}f_{h-2}(S_{h-3})f_{h-1}(I_{h-2}),\label{eq:h-j2}\end{equation}
for $h=1,2,\dots,k$. From (\ref{eq:h-j2}) we obtain\begin{align*}
I_{h} & \geq I_{h-2}\left(f_{h-2}(S_{h-3})f_{h-1}(I_{h-2})\right)\\
 & \geq I_{h-4}\left(f_{h-4}(S_{h-5})f_{h-3}(I_{h-4})\right)\left(f_{h-2}(S_{h-3})f_{h-1}(I_{h-2})\right)\\
 & \geq\dots\\
 & \geq I_{h-k}\underline{G}_{h},\end{align*}
 where \begin{align}
\underline{G}_{h} & =\left(f_{h-k}(S_{h-k-1})f_{h+1-k}(I_{h-k})\right)\dots\left(f_{h-4}(S_{h-5})f_{h-3}(I_{h-4})\right)\left(f_{h-2}(S_{h-3})f_{h-1}(I_{h-2})\right).\label{eq:underg}\end{align}
Then \[
1\geq\underline{G}_{h},\]
holds for $h=1,2,\dots,k$, because we have $I_{h}=I_{h-k}$ from
(\ref{periodicityofsi}). Consequently, it holds that\begin{equation}
\underline{G}_{h}\leq1\leq\overline{G}_{h},\text{ for }h=1,2,\dots,k.\label{eq:Grel1}\end{equation}
By (\ref{eq:overg}) and (\ref{eq:underg}), we see\begin{equation}
\overline{G}_{1}=\underline{G}_{k}\text{ and }\overline{G}_{h}=\underline{G}_{h-1}\text{ for }h\in\left\{ 2,3,\dots,k\right\} .\label{eq:Grel2}\end{equation}
Therefore, we obtain (\ref{eq:concl_lem_even}) from (\ref{eq:Grel1})
and (\ref{eq:Grel2}). The proof is complete.\qed 
\end{pf}
We show that there exists a $k$-periodic solution which is globally
attractive for any solution for the case where $k$ is an even integer.
\begin{thm}
\label{thm:Theore,_even}Let $k$ be an even integer. Assume that
(\ref{eq:ass00}) and (\ref{eq:assumption}). If (\ref{eq:theorem11}),
then there exists a $k$-periodic solution $x_{n}^{*}$ such that
$x_{n}^{*}=x_{n+k}^{*}$ which is globally attractive, that is, for
any solution of (\ref{eq:Mother_generalPielou}), it holds that\[
\lim_{n\to+\infty}\left(x_{n}-x_{n}^{*}\right)=0.\]
\end{thm}
\begin{pf}
To obtain the conclusion, we will show \begin{equation}
S_{h}=I_{h}\text{ for }h=1,2,\dots,k.\label{eq:concl_even}\end{equation}
Let\[
U_{h,h-j}=\lim_{m\to+\infty}x_{k\overline{n}_{m}^{h}+h-j}\text{ for }j=2,3\dots,k+1,\]
and we claim that it holds

\begin{equation}
U_{h,h-j}=\begin{cases}
S_{h-j} & \text{ for }j=2,4,\dots,k,\\
I_{h-j} & \text{ for }j=3,5,\dots,k+1,\end{cases}\label{inductive}\end{equation}
for $h=1,2,\dots,k$.

From (\ref{eq:Mother_generalPielou}), it holds\begin{align}
x_{k\overline{n}_{m}^{h}+h-j} & =x_{k\overline{n}_{m}^{h}+h-j-2}f_{k\overline{n}_{m}^{h}+h-j-2}(x_{k\overline{n}_{m}^{h}+h-j-3})f_{k\overline{n}_{m}^{h}+h-j-1}(x_{k\overline{n}_{m}^{h}+h-j-2})\nonumber \\
 & =x_{k\overline{n}_{m}^{h}+h-j-2}f_{h-j-2}(x_{k\overline{n}_{m}^{h}+h-j-3})f_{h-j-1}(x_{k\overline{n}_{m}^{h}+h-j-2}).\label{eq:gp2_even_gen}\end{align}

Firstly, we show \begin{equation}
U_{h,h-j}=\begin{cases}
S_{h-j} & \text{ for }j=2,\\
I_{h-j} & \text{ for }j=3,\end{cases}\label{eq:U1_gen1}\end{equation}
for $h=1,2,\dots,k$. Suppose that there exists $\overline{h}$ such
that $U_{\overline{h},\overline{h}-2}<S_{\overline{h}-2},$ or $U_{\overline{h},\overline{h}-3}>I_{\overline{h}-3}$.
By letting $m\to+\infty$ and considering the limiting equation of
(\ref{eq:gp2_even_gen}) with $j=0$, we obtain\begin{equation}
S_{h}=U_{h,h-2}f_{h-2}(U_{h,h-3})f_{h-1}(U_{h,h-2}),\label{eq:Sh1_gen}\end{equation}
for $h=1,2,\dots,k$. Then, we see that one of the following holds\[
\Pi_{i=1}^{p}S_{2i}<\Pi_{i=1}^{p}\left[S_{2i-2}f_{2i-2}(I_{2i-3})f_{2i-1}(S_{2i-2})\right],\]
or\[
\Pi_{i=0}^{p-1}S_{2i+1}<\Pi_{i=0}^{p-1}\left[S_{2i-1}f_{2i-1}(I_{2i-2})f_{2i}(S_{2i-1})\right],\]
where $p=\frac{k}{2}$, since \[
\begin{cases}
\Pi_{i=1}^{p}S_{2i} & =\Pi_{i=1}^{p}\left[U_{2i,2i-2}f_{2i-2}(U_{2i,2i-3})f_{2i-1}(U_{2i,2i-2})\right],\\
\Pi_{i=0}^{p-1}S_{2i+1} & =\Pi_{i=0}^{p-1}\left[U_{2i+1,2i-1}f_{2i-1}(U_{2i+1,2i-2})f_{2i}(U_{2i+1,2i-1})\right],\end{cases}\]
follows from (\ref{eq:Sh1_gen}). Then, we obtain \[
1<\Pi_{i=1}^{p}\left[f_{2i-2}(I_{2i-3})f_{2i-1}(S_{2i-2})\right],\]
or \[
1<\Pi_{i=0}^{p-1}\left[f_{2i-1}(I_{2i-2})f_{2i}(S_{2i-1})\right],\]
and this gives a contradiction to (\ref{eq:concl_lem_even}) with
$h=k$ and $h=k-1$, respectively, in Lemma \ref{lem:Lemma_even}.
Thus (\ref{eq:U1_gen1}) holds.

Next, we assume that \begin{equation}
U_{h,h-j}=\begin{cases}
S_{h-j} & \text{ for }j=j_{1},\\
I_{h-j} & \text{ for }j=j_{1}+1,\end{cases}\label{eq:U1_gen2}\end{equation}
for $h=1,2,\dots,k$ where $j_{1}$ is a positive even number. Under
the assumption (\ref{eq:U1_gen2}), we show \begin{equation}
U_{h,h-j}=\begin{cases}
S_{h-j} & \text{ for }j=j_{1}+2,\\
I_{h-j} & \text{ for }j=j_{1}+3,\end{cases}\label{eq:U1_gen}\end{equation}
for $h=1,2,\dots,k$. Suppose that there exists $\overline{h}$ such
that $U_{\overline{h},\overline{h}-j_{1}-2}<S_{\overline{h}-j_{1}-2},$
or $U_{\overline{h},\overline{h}-j_{1}-3}>I_{\overline{h}-j_{1}-3}$.
By considering the limiting equation of (\ref{eq:gp2_even_gen}) with
$j=j_{1}$ and substituting (\ref{eq:U1_gen2}),\[
U_{h,h-j_{1}}=S_{h-j_{1}}=U_{h,h-j_{1}-2}f_{h-j_{1}-2}(U_{h,h-j_{1}-3})f_{h-j_{1}-1}(U_{h,h-j_{1}-2}),\]
for $h=1,2,\dots,k$. Then, similar to the above discussion, we can
show that (\ref{eq:U1_gen}) holds. Thus, (\ref{inductive}) holds
by mathematical induction.

From (\ref{eq:Mother_generalPielou}), it holds \begin{align*}
x_{k\overline{n}_{m}^{h}+h} & =x_{k\overline{n}_{m}^{h}+h-(k-1)}f_{k\overline{n}_{m}^{h}+h-(k-1)}(x_{k\overline{n}_{m}^{h}+h-k})\dots f_{k\overline{n}_{m}^{h}+h-2}(x_{k\overline{n}_{m}^{h}+h-3})f_{k\overline{n}_{m}^{h}+h-1}(x_{k\overline{n}_{m}^{h}+h-2})\\
 & =x_{k\overline{n}_{m}^{h}+h-(k-1)}f_{h-(k-1)}(x_{k\overline{n}_{m}^{h}+h-k})\dots f_{h-2}(x_{k\overline{n}_{m}^{h}+h-3})f_{h-1}(x_{k\overline{n}_{m}^{h}+h-2}),\end{align*}
and by considering the limiting equation and using (\ref{inductive}),
we obtain \begin{align*}
S_{h} & =U_{h,h-(k-1)}f_{h-(k-1)}(U_{h,h-k})\dots f_{h-2}(U_{h,h-3})f_{h-1}(U_{h,h-2})\\
 & =I_{h-(k-1)}f_{h-(k-1)}(S_{h-k})\dots f_{h-2}(I_{h-3})f_{h-1}(S_{h-2}).\end{align*}
By (\ref{eq:concl_lem_even}) in Lemma \ref{lem:Lemma_even}, we see
\[
f_{h-(k-1)}(S_{h-k})\dots f_{h-2}(I_{h-3})f_{h-1}(S_{h-2})=\frac{1}{f_{h-k}(I_{h-k-1})},\]
hence, it holds that \begin{equation}
S_{h}f_{h-k}(I_{h-k-1})=I_{h+1},\label{Sendproof2p}\end{equation}
for $h=1,2,\dots,k$. Similar to the above discussion, it also holds
\begin{equation}
I_{h}f_{h-k}(S_{h-k-1})=S_{h+1},\label{Iendproof2p}\end{equation}
for $h=1,2,\dots,k$. Consequently, it holds \[
I_{h+1}=S_{h}f_{h-k}(I_{h-k-1})\leq I_{h}f_{h-k}(S_{h-k-1})=S_{h+1},\]
by (\ref{Sendproof2p}) and (\ref{Iendproof2p}), and hence, (\ref{eq:concl_even})
holds. Then, from (\ref{eq:defsh}), we see \[
\liminf_{n\to+\infty}x_{kn+h}=\limsup_{n\to+\infty}x_{kn+h}=\lim_{n\to+\infty}x_{kn+h},\]
for $h=1,2,\dots,k$, and there exist $k$ positive constants $x_{h}^{*}=S_{h}=I_{h},h=1,2,\dots,k$
such that \[
x_{h}^{*}=\lim_{n\to+\infty}x_{kn+h}.\]
The proof is complete.\qed 
\end{pf}

\section{Global attractivity for the case where $k$ is an arbitrary odd integer}

In this section, we show that Theorem \ref{thm:mainthm} holds when
$k$ is an odd integer. For the reader, we first consider the case
$k=3$ in Section 5.1. (\ref{eq:concl_lem_3}) in Lemma \ref{lem:Lemma_3}
has an important role in this subsection. Then, we give Theorem \ref{thm:Theorem_3}
which states that there exists a $3$-periodic solution which is globally
attractive. In Section 5.2, we generalize these results to the case
where $k$ is an arbitrary odd integer.

\subsection{Case: $k=3$}

First, we introduce the following lemma which plays a crucial role
in the proof of Theorem \ref{thm:Theorem_3}.
\begin{lem}
\label{lem:Lemma_3}Let $k=3$. Assume that (\ref{eq:ass00}) and
(\ref{eq:assumption}). If (\ref{eq:theorem11}), then it holds that\begin{equation}
\begin{cases}
S_{h} & =S_{h-2}f_{h-2}(I_{h-3})f_{h-1}(S_{h-2}),\\
I_{h} & =I_{h-2}f_{h-2}(S_{h-3})f_{h-1}(I_{h-2}),\end{cases}\label{eq:concl_lem_3}\end{equation}
for $h=1,2,3$.\end{lem}
\begin{pf}
From (\ref{eq:Mother_generalPielou}), it holds\begin{align*}
x_{kn+h} & =x_{kn+h-2}f_{kn+h-2}(x_{kn+h-3})f_{kn+h-1}(x_{kn+h-2})\\
 & =x_{kn+h-2}f_{h-2}(x_{kn+h-3})f_{h-1}(x_{kn+h-2}),\end{align*}
for $h=1,2,3$. Then it follows that\[
\begin{cases}
x_{k\overline{n}_{m}^{h}+h} & =x_{k\overline{n}_{m}^{h}+h-2}f_{h-2}(x_{k\overline{n}_{m}^{h}+h-3})f_{h-1}(x_{k\overline{n}_{m}^{h}+h-2}),\\
x_{k\underline{n}_{m}^{h}+h} & =x_{k\underline{n}_{m}^{h}+h-2}f_{h-2}(x_{k\underline{n}_{m}^{h}+h-3})f_{h-1}(x_{k\underline{n}_{m}^{h}+h-2}),\end{cases}\]
for $h=1,2,3$. By considering the limiting equation and using (\ref{eq:ass00})
and (\ref{eq:assumption}), we get\begin{equation}
\begin{cases}
S_{h} & \leq S_{h-2}f_{h-2}(I_{h-3})f_{h-1}(S_{h-2}),\\
I_{h} & \geq I_{h-2}f_{h-2}(S_{h-3})f_{h-1}(I_{h-2}),\end{cases}\label{3lemma}\end{equation}
for $h=1,2,3$. In order to obtain the conclusion, we show that (\ref{3lemma})
holds with equality.

From (\ref{3lemma}), we see that \begin{equation}
\begin{cases}
\frac{S_{h}f_{h-1}(I_{h-2})}{I_{h}f_{h-1}(S_{h-2})} & \leq\frac{S_{h-2}}{I_{h}}f_{h-2}(I_{h-3})f_{h-1}(I_{h-2}),\\
\frac{S_{h-2}f_{h-2}(I_{h-3})}{I_{h-2}f_{h-2}(S_{h-3})} & \geq\frac{S_{h-2}}{I_{h}}f_{h-2}(I_{h-3})f_{h-1}(I_{h-2}),\end{cases}\label{3lemma4}\end{equation}
for $h=1,2,3$. It then follows \begin{equation}
\frac{S_{h}f_{h-1}(I_{h-2})}{I_{h}f_{h-1}(S_{h-2})}\leq\frac{S_{h-2}f_{h-2}(I_{h-3})}{I_{h-2}f_{h-2}(S_{h-3})},\label{3lemma5}\end{equation}
 for $h=1,2,3$. By multiplying (\ref{3lemma5}), we obtain \[
\left(\frac{S_{h}f_{h-1}(I_{h-2})}{I_{h}f_{h-1}(S_{h-2})}\right)\left(\frac{S_{h+1}f_{h}(I_{h-1})}{I_{h+1}f_{h}(S_{h-1})}\right)\leq\left(\frac{S_{h-2}f_{h-2}(I_{h-3})}{I_{h-2}f_{h-2}(S_{h-3})}\right)\left(\frac{S_{h-1}f_{h-1}(I_{h-2})}{I_{h-1}f_{h-1}(S_{h-2})}\right),\]
and hence, it holds that\[
\frac{S_{h}S_{h+1}f_{h}(I_{h-1})}{I_{h}I_{h+1}f_{h}(S_{h-1})}\leq\frac{S_{h-2}S_{h-1}f_{h-2}(I_{h-3})}{I_{h-2}I_{h-1}f_{h-2}(S_{h-3})},\]
for $h=1,2,3$. Therefore, the following inequalities hold\begin{align*}
\begin{cases}
\frac{S_{1}S_{2}f_{1}(I_{3})}{I_{1}I_{2}f_{1}(S_{3})} & \leq\frac{S_{2}S_{3}f_{2}(I_{1})}{I_{2}I_{3}f_{2}(S_{1})},\\
\frac{S_{2}S_{3}f_{2}(I_{1})}{I_{2}I_{3}f_{2}(S_{1})} & \leq\frac{S_{3}S_{1}f_{3}(I_{2})}{I_{3}I_{1}f_{3}(S_{2})},\\
\frac{S_{3}S_{1}f_{3}(I_{2})}{I_{3}I_{1}f_{3}(S_{2})} & \leq\frac{S_{1}S_{2}f_{1}(I_{3})}{I_{1}I_{2}f_{1}(S_{3})},\end{cases}\end{align*}
and it follows\[
\frac{S_{1}S_{2}f_{1}(I_{3})}{I_{1}I_{2}f_{1}(S_{3})}\leq\frac{S_{2}S_{3}f_{2}(I_{1})}{I_{2}I_{3}f_{2}(S_{1})}\leq\frac{S_{3}S_{1}f_{3}(I_{2})}{I_{3}I_{1}f_{3}(S_{2})}\leq\frac{S_{1}S_{2}f_{1}(I_{3})}{I_{1}I_{2}f_{1}(S_{3})}.\]
Hence, we see that \begin{equation}
\frac{S_{1}S_{2}f_{1}(I_{3})}{I_{1}I_{2}f_{1}(S_{3})}=\frac{S_{2}S_{3}f_{2}(I_{1})}{I_{2}I_{3}f_{2}(S_{1})}=\frac{S_{3}S_{1}f_{3}(I_{2})}{I_{3}I_{1}f_{3}(S_{2})}=\frac{S_{1}S_{2}f_{1}(I_{3})}{I_{1}I_{2}f_{1}(S_{3})}.\label{3lemma7}\end{equation}
(\ref{3lemma7}) implies that (\ref{3lemma4}) holds with equality
and thus, (\ref{3lemma}) also holds with equality for $h=1,2,3$.
Therefore, we obtain the conclusion. The proof is complete.\qed
\end{pf}
Then, we obtain the following result.

\begin{thm}
\label{thm:Theorem_3}Let $k=3$. Assume that (\ref{eq:ass00}) and
(\ref{eq:assumption}). If (\ref{eq:theorem11}), then there exists
a $3$-periodic solution $x_{n}^{*}$ such that $x_{n}^{*}=x_{n+3}^{*}$
which is globally attractive, that is, for any solution of (\ref{eq:Mother_generalPielou}),
it holds that\[
\lim_{n\to+\infty}\left(x_{n}-x_{n}^{*}\right)=0.\]
\end{thm}
\begin{pf}
In order to obtain the conclusion, we will show \begin{equation}
S_{h}=I_{h},h=1,2,3.\label{eq:eq:concl_3}\end{equation}
Let\[
U_{k,k-j}=\lim_{m\to+\infty}x_{k\overline{n}_{m}^{k}+k-j}\text{ for }j=1,2,\dots,5,\]
and we claim\begin{equation}
U_{k,k-j}=\begin{cases}
S_{k-j} & \text{ for }j=2,4,\\
I_{k-j} & \text{ for }j=1,3,5.\end{cases}\label{eq:claim3_0}\end{equation}

At first, we see that it holds\begin{align}
x_{k\overline{n}_{m}^{k}+k-j} & =x_{k\overline{n}_{m}^{k}+k-j-2}f_{k\overline{n}_{m}^{k}+k-j-2}(x_{k\overline{n}_{m}^{k}+k-j-3})f_{k\overline{n}_{m}^{k}+k-j-1}(x_{k\overline{n}_{m}^{k}+k-j-2})\nonumber \\
 & =x_{k\overline{n}_{m}^{k}+k-j-2}f_{k-j-2}(x_{k\overline{n}_{m}^{k}+k-j-3})f_{k-j-1}(x_{k\overline{n}_{m}^{k}+k-j-2}),\label{eq:ite3}\end{align}
for $j=0,1,2,\dots,$ from (\ref{eq:Mother_generalPielou}).

Firstly, we show\begin{equation}
U_{k,k-j}=\begin{cases}
S_{k-j} & \text{ for }j=2,\\
I_{k-j} & \text{ for }j=3.\end{cases}\label{eq:claim3_1}\end{equation}
Suppose that $U_{k,k-2}<S_{k-2}$ or $U_{k,k-3}>I_{k-3}$. By considering
the limiting equation of (\ref{eq:ite3}) with $j=0$, it follows\[
S_{k}=U_{k,k-2}f_{k-2}(U_{k,k-3})f_{k-1}(U_{k,k-2}).\]
Then, by (\ref{eq:ass00}) and (\ref{eq:assumption}), it follows
\[
S_{k}<S_{k-2}f_{k-2}(I_{k-3})f_{k-1}(S_{k-2}).\]
This gives a contradiction to (\ref{eq:concl_lem_3}) with $h=k$
in Lemma \ref{lem:Lemma_3}. Thus, (\ref{eq:claim3_1}) holds. 

Next, we show \begin{equation}
U_{k,k-j}=\begin{cases}
S_{k-j} & \text{ for }j=4,\\
I_{k-j} & \text{ for }j=5.\end{cases}\label{eq:claim3_2}\end{equation}
Suppose that $U_{k,k-4}<S_{k-4}$ or $U_{k,k-5}>I_{k-5}$. By considering
the limiting equation of (\ref{eq:ite3}) with $j=2$ and substituting
(\ref{eq:claim3_1}), it follows\[
U_{k,k-2}=S_{k-2}=U_{k,k-4}f_{k-4}(U_{k,k-5})f_{k-3}(U_{k,k-4}).\]
Then, by (\ref{eq:ass00}) and (\ref{eq:assumption}), it follows
\[
S_{k-2}<S_{k-4}f_{k-4}(I_{k-5})f_{k-3}(S_{k-4}).\]
This gives a contradiction to (\ref{eq:concl_lem_3}) with $h=k-2$
in Lemma \ref{lem:Lemma_3}. Thus, (\ref{eq:claim3_2}) holds. 

By considering the limiting equation of (\ref{eq:ite3}) with $j=1$
and using (\ref{eq:claim3_1})-(\ref{eq:claim3_2}), it follows \begin{align*}
U_{k,k-1} & =U_{k,k-3}f_{k-3}(U_{k,k-4})f_{k-2}(U_{k,k-3})\\
 & =I_{k-3}f_{k-3}(S_{k-4})f_{k-2}(I_{k-3}).\end{align*}
Hence, it holds \begin{equation}
U_{k,k-1}=I_{k-1},\label{eq:claim3_3}\end{equation}
by (\ref{eq:concl_lem_3}) in Lemma \ref{lem:Lemma_3}. Consequently,
(\ref{eq:claim3_0}) holds from (\ref{eq:claim3_1}), (\ref{eq:claim3_2})
and (\ref{eq:claim3_3}).

From (\ref{eq:Mother_generalPielou}), it holds\[
x_{k\overline{n}_{m}^{k}+k}=x_{k\overline{n}_{m}^{k}+k-1}f_{k\overline{n}_{m}^{k}+k-1}(x_{k\overline{n}_{m}^{k}+k-2})=x_{k\overline{n}_{m}^{k}+k-1}f_{k-1}(x_{k\overline{n}_{m}^{k}+k-2}),\]
 and by considering the limiting equation and using (\ref{eq:claim3_0}),
we obtain \begin{equation}
S_{k}=I_{k-1}f_{k-1}(S_{k-2}).\label{s1}\end{equation}
Similar to the above discussion, it also holds\begin{equation}
I_{k}=S_{k-1}f_{k-1}(I_{k-2}).\label{i1}\end{equation}
Consequently, it holds \[
I_{k}=S_{k-1}f_{k-1}(I_{k-2})\leq I_{k-1}f_{k-1}(S_{k-2})=S_{k},\]
by (\ref{s1}) and (\ref{i1}), and hence (\ref{eq:eq:concl_3}) holds.
Then, \[
\liminf_{n\to+\infty}x_{kn+h}=\limsup_{n\to+\infty}x_{kn+h}=\lim_{n\to+\infty}x_{kn+h},\]
for $h=1,2,3$ and there exist $3$-positive constants $x_{h}^{*}=S_{h}=I_{h},h=1,2,3$
such that \[
x_{h}^{*}=\lim_{n\to+\infty}x_{kn+h}.\]
The proof is complete.\qed 
\end{pf}

\subsection{Case: $k$ is an odd integer}
\begin{lem}
\label{lem:Lemma_odd} Let $k$ be an odd integer. Assume that (\ref{eq:ass00})
and (\ref{eq:assumption}). If (\ref{eq:theorem11}), then it holds
that\begin{equation}
\begin{cases}
S_{h} & =S_{h-2}f_{h-2}(I_{h-3})f_{h-1}(S_{h-2}),\\
I_{h} & =I_{h-2}f_{h-2}(S_{h-3})f_{h-1}(I_{h-2}),\end{cases}\label{eq:concl_lem_odd}\end{equation}
for $h=1,2,\dots,k$.\end{lem}
\begin{pf}
From (\ref{eq:Mother_generalPielou}), it holds \begin{align*}
x_{kn+h} & =x_{kn+h-2}f_{kn+h-2}(x_{kn+h-3})f_{kn+h-1}(x_{kn+h-2})\\
 & =x_{kn+h-2}f_{h-2}(x_{kn+h-3})f_{h-1}(x_{kn+h-2}),\end{align*}
for $h=1.2.\dots,k$. Then, it follows that\[
\begin{cases}
x_{k\overline{n}_{m}^{h}+h} & =x_{k\overline{n}_{m}^{h}+h-2}f_{h-2}(x_{k\overline{n}_{m}^{h}+h-3})f_{h-1}(x_{k\overline{n}_{m}^{h}+h-2}),\\
x_{k\underline{n}_{m}^{h}+h} & =x_{k\underline{n}_{m}^{h}+h-2}f_{h-2}(x_{k\underline{n}_{m}^{h}+h-3})f_{h-1}(x_{k\underline{n}_{m}^{h}+h-2}),\end{cases}\]
for $h=1,2,\dots,k$. By considering the limiting equation and using
(\ref{eq:ass00}) and (\ref{eq:assumption}), we get\begin{equation}
\begin{cases}
S_{h} & \leq S_{h-2}f_{h-2}(I_{h-3})f_{h-1}(S_{h-2}),\\
I_{h} & \geq I_{h-2}f_{h-2}(S_{h-3})f_{h-1}(I_{h-2}),\end{cases}\label{3@lemma}\end{equation}
for $h=1,2,\dots,k$. In order to obtain the conclusion, we show that
(\ref{3@lemma}) holds with equality.

From (\ref{3@lemma}), we see that \begin{equation}
\begin{cases}
\frac{S_{h}f_{h-1}(I_{h-2})}{I_{h}f_{h-1}(S_{h-2})} & \leq\frac{S_{h-2}}{I_{h}}f_{h-2}(I_{h-3})f_{h-1}(I_{h-2}),\\
\frac{S_{h-2}f_{h-2}(I_{h-3})}{I_{h-2}f_{h-2}(S_{h-3})} & \geq\frac{S_{h-2}}{I_{h}}f_{h-2}(I_{h-3})f_{h-1}(I_{h-2}),\end{cases}\label{3@lemma4}\end{equation}
for $h=1,2,\dots,k$. Thus, it follows that \begin{equation}
\frac{S_{h}f_{h-1}(I_{h-2})}{I_{h}f_{h-1}(S_{h-2})}\leq\frac{S_{h-2}f_{h-2}(I_{h-3})}{I_{h-2}f_{h-2}(S_{h-3})},\label{3@lemma5}\end{equation}
for $h=1,2,\dots,k$. By multiplying (\ref{3@lemma5}), we obtain
\[
\left(\frac{S_{h}f_{h-1}(I_{h-2})}{I_{h}f_{h-1}(S_{h-2})}\right)\left(\frac{S_{h+1}f_{h}(I_{h-1})}{I_{h+1}f_{h}(S_{h-1})}\right)\leq\left(\frac{S_{h-2}f_{h-2}(I_{h-3})}{I_{h-2}f_{h-2}(S_{h-3})}\right)\left(\frac{S_{h-1}f_{h-1}(I_{h-2})}{I_{h-1}f_{h-1}(S_{h-2})}\right),\]
and hence, it holds that\[
\frac{S_{h}S_{h+1}f_{h}(I_{h-1})}{I_{h}I_{h+1}f_{h}(S_{h-1})}\leq\frac{S_{h-2}S_{h-1}f_{h-2}(I_{h-3})}{I_{h-2}I_{h-1}f_{h-2}(S_{h-3})},\]
for $h=1,2,\dots,k$. Therefore, the following inequalities hold\[
\begin{cases}
\frac{S_{1}S_{2}f_{1}(I_{k})}{I_{1}I_{2}f_{1}(S_{k})} & \leq\frac{S_{k-1}S_{k}f_{k-1}(I_{k-2})}{I_{k-1}I_{k}f_{k-1}(S_{k-2})},\\
\frac{S_{k-1}S_{k}f_{k-1}(I_{k-2})}{I_{k-1}I_{k}f_{k-1}(S_{k-2})} & \leq\frac{S_{k-3}S_{k-2}f_{k-3}(I_{k-4})}{I_{k-3}I_{k-2}f_{k-3}(S_{k-4})},\\
\dots & \dots\\
\frac{S_{k}S_{k+1}f_{k}(I_{k-1})}{I_{k}I_{k+1}f_{k}(S_{k-1})} & \leq\frac{S_{k-2}S_{k-1}f_{k-2}(I_{k-3})}{I_{k-2}I_{k-1}f_{k-2}(S_{k-3})},\\
\dots & \dots\\
\frac{S_{3}S_{4}f_{3}(I_{2})}{I_{3}I_{4}f_{3}(S_{2})} & \leq\frac{S_{1}S_{2}f_{1}(I_{k})}{I_{1}I_{2}f_{1}(S_{k})}.\end{cases}\]
Thus, we obtain \begin{equation}
\frac{S_{1}S_{2}f_{1}(I_{k})}{I_{1}I_{2}f_{1}(S_{k})}=\dots=\frac{S_{k-1}S_{k}f_{k-1}(I_{k-2})}{I_{k-1}I_{k}f_{k-1}(S_{k-2})}=\frac{S_{k}S_{k+1}f_{k}(I_{k-1})}{I_{k}I_{k+1}f_{k}(S_{k-1})}.\label{3@lemma7}\end{equation}
(\ref{3@lemma7}) implies that (\ref{3@lemma4}) holds with equality
and thus, (\ref{3@lemma}) also holds with equality for $h=1,2,\dots,k$.
Then, we obtain the conclusion. Hence, the proof is complete.\qed\end{pf}
\begin{thm}
\label{thm:Theorem_odd}Let $k$ be an odd integer. Assume that (\ref{eq:ass00})
and (\ref{eq:assumption}). If (\ref{eq:theorem11}), then there exists
a $k$-periodic solution $x_{n}^{*}$ such that $x_{n}^{*}=x_{n+k}^{*}$
which is globally attractive, that is, for any solution of (\ref{eq:Mother_generalPielou}),
it holds that\[
\lim_{n\to+\infty}\left(x_{n}-x_{n}^{*}\right)=0.\]
\end{thm}
\begin{pf}
In order to obtain the conclusion, we will show \begin{equation}
S_{h}=I_{h},h=1,2,\dots,k.\label{eq:concl_odd}\end{equation}
 Let\[
U_{k,k-j}=\lim_{m\to+\infty}x_{k\overline{n}_{m}^{k}+k-j}\text{ for }j=1,2,\dots.,k,\]
and we claim \begin{equation}
U_{k,k-j}=\begin{cases}
S_{k-j}, & \text{ for }j=2,4,\dots,k-1,\\
I_{k-j}, & \text{ for }j=1,3,5,\dots,k.\end{cases}\label{inductive2p+1}\end{equation}

At first, we see that it holds\begin{align}
x_{k\overline{n}_{m}^{k}+k-j} & =x_{k\overline{n}_{m}^{k}+k-j-2}f_{k\overline{n}_{m}^{k}+k-j-2}(x_{k\overline{n}_{m}^{k}+k-j-3})f_{k\overline{n}_{m}^{k}+k-j-1}(x_{k\overline{n}_{m}^{k}+k-j-2})\nonumber \\
 & =x_{k\overline{n}_{m}^{k}+k-j-2}f_{k-j-2}(x_{k\overline{n}_{m}^{k}+k-j-3})f_{k-j-1}(x_{k\overline{n}_{m}^{k}+k-j-2}),\label{eq:ite_odd}\end{align}
for $j=0,1,2,\dots$.

Firstly, we show \begin{equation}
U_{k,k-j}=\begin{cases}
S_{k-j} & \text{ for }j=2,\\
I_{k-j} & \text{ for }j=3.\end{cases}\label{inductive1}\end{equation}
Suppose that $U_{k,k-2}<S_{k-2}$ or $U_{k,k-3}>I_{k-3}$. By considering
the limiting equation of (\ref{eq:ite_odd}) with $j=0$, it follows\[
S_{k}=U_{k,k-2}f_{k-2}(U_{k,k-3})f_{k-1}(U_{k,k-2}).\]
Then, by (\ref{eq:ass00}) and (\ref{eq:assumption}), it follows
\[
S_{k}<S_{k-2}f_{k-2}(I_{k-3})f_{k-1}(S_{k-2}).\]
This gives a contradiction to (\ref{eq:concl_lem_odd}) with $h=k$
in Lemma \ref{lem:Lemma_odd}. Thus, (\ref{inductive1}) holds. 

Next we assume that\begin{equation}
U_{k,k-j}=\begin{cases}
S_{k-j} & \text{ for }j=j_{1},\\
I_{k-j} & \text{ for }j=j_{1}+1,\end{cases}\label{inductive2}\end{equation}
where $j_{1}$ is a positive even integer. Under the assumption (\ref{inductive2}),
we show\begin{equation}
U_{k,k-j}=\begin{cases}
S_{k-j} & \text{ for }j=j_{1}+2,\\
I_{k-j} & \text{ for }j=j_{1}+3.\end{cases}\label{inductive3}\end{equation}
Suppose that $U_{k,k-j_{1}-2}<S_{k-j_{1}-2}$ or $U_{k,k-j_{1}-3}>I_{k-j_{1}-3}$.
By considering the limiting equation of (\ref{eq:ite_odd}) with $j=j_{1}$,
it follows\[
U_{k,k-j_{1}}=S_{k-j_{1}}=U_{k,k-j_{1}-2}f_{k-j_{1}-2}(U_{k,k-j_{1}-3})f_{k-j_{1}-1}(U_{k,k-j_{1}-2}).\]
Then, by (\ref{eq:ass00}) and (\ref{eq:assumption}), it follows
\[
S_{k-j_{1}}<S_{k-j_{1}-2}f_{k-j_{1}-2}(I_{k-j_{1}-3})f_{k-1}(S_{k-j_{1}-2}).\]
This gives a contradiction to (\ref{eq:concl_lem_odd}) with $h=k-j_{1}$
in Lemma \ref{lem:Lemma_odd}. Thus, (\ref{inductive3}) holds. 

By considering the limiting equation of (\ref{eq:ite_odd}) with $j=1$
and substituting (\ref{inductive1}) and (\ref{inductive3}) with
$j_{1}=2$, it follows \begin{align*}
U{}_{k,k-1} & =U_{k,k-3}f_{k-3}(U_{k,k-4})f_{k-2}(U_{k,k-3})\\
 & =I_{k-3}f_{k-3}(S_{k-4})f_{k-2}(I_{k-3}).\end{align*}
Hence, \begin{equation}
U_{k,k-1}=I_{k-1},\label{eq:inductive4}\end{equation}
by (\ref{eq:concl_lem_odd}) with $h=k-1$ in Lemma \ref{lem:Lemma_odd}.
Hence, (\ref{inductive2p+1}) holds by (\ref{eq:inductive4}) and
mathematical induction.

From (\ref{eq:Mother_generalPielou}), it holds\[
x_{k\overline{n}_{m}^{k}+k-j}=x_{k\overline{n}_{m}^{k}+k-j-1}f_{k\overline{n}_{m}^{k}+k-j-1}(x_{k\overline{n}_{m}^{k}+k-j-2})=x_{k\overline{n}_{m}^{k}+k-j-1}f_{k-j-1}(x_{k\overline{n}_{m}^{k}+k-j-2}),\]
for $j=0,1,2,\dots,k-1$, and by considering the limiting equation
and using (\ref{inductive2p+1}), we obtain \begin{equation}
\begin{cases}
I_{k-j} & =S_{k-j-1}f_{k-j-1}(I_{k-j-2}),j=1,3,\dots,k,\\
S_{k-j} & =I_{k-j-1}f_{k-j-1}(S_{k-j-2}),j=2,4,\dots,k-1.\end{cases}\label{s}\end{equation}

Similar to the above discussion, it also holds\begin{equation}
\begin{cases}
S_{k-j} & =I_{k-j-1}f_{k-j-1}(S_{k-j-2}),j=1,3,\dots,k,\\
I_{k-j} & =S_{k-j-1}f_{k-j-1}(I_{k-j-2}),j=2,4,\dots,k-1.\end{cases}\label{i}\end{equation}
Consequently, it holds \[
I_{k-j}=S_{k-j-1}f_{k-j-1}(I_{k-j-2})\leq I_{k-j-1}f_{k-j-1}(S_{k-j-2})=S_{k-j},\]
for $j=1,2,\dots,k$, by (\ref{s}) and (\ref{i}), and hence (\ref{eq:concl_odd})
holds. Then, \[
\liminf_{n\to+\infty}x_{kn+h}=\limsup_{n\to+\infty}x_{kn+h}=\lim_{n\to+\infty}x_{kn+h},\]
for $h=1,2,\dots,k$ and there exist $k$ positive constants $x_{h}^{*}=S_{h}=I_{h},h=1,2,\dots,k$
such that \[
x_{h}^{*}=\lim_{n\to+\infty}x_{kn+h}.\]
The proof is complete.\qed 

Finally, we establish Theorem \ref{thm:mainthm} from Theorems \ref{thm:Theorem_2},
\ref{thm:Theore,_even}, \ref{thm:Theorem_3} and \ref{thm:Theorem_odd}.
\end{pf}

\section{Applications}

In this section, we give two examples to demonstrate our result. At
first, we introduce the following example.

\begin{equation}
f_{n}(x)=\frac{\lambda}{1+(\lambda-1)\frac{x}{K_{n}}},n=0,1,2\dots,x\in[0,+\infty),\label{eq:beverton}\end{equation}
where $\lambda>1$ and $K_{n}=K_{n+k}>0$ for $n=0,1,2,\dots,$ ($k$
is an arbitrary positive integer). (\ref{eq:Mother_generalPielou})
with (\ref{eq:beverton}) is called a delayed Beverton-Holt equation
(see also \cite{MR2397325,MR2151684}). Since $f_{n}(x)$ is periodic
on $n$ and satisfies (\ref{eq:ass00}) and (\ref{eq:assumption}),
the following result is derived from Theorem \ref{thm:mainthm}.
\begin{cor}
\label{cor:bev}There exists a $k$-periodic solution $x_{n}^{*}$
such that $x_{n}^{*}=x_{n+k}^{*}$ which is globally attractive, that
is, for any solution of (\ref{eq:Mother_generalPielou}) with (\ref{eq:beverton}),
it holds that \[
\lim_{n\to+\infty}(x_{n}-x_{n}^{*})=0.\]

\end{cor}
One can see (\ref{eq:Mother_generalPielou}) with (\ref{eq:beverton})
is equivalent to the Pielou's equation by a simple transformation.
Therefore, Corollary \ref{cor:bev} is also derived from Theorem A.

Secondly, we introduce the following example.\begin{equation}
f_{n}(x)=\frac{\beta_{n}}{1+\alpha_{n}^{1}\frac{x}{1+\alpha_{n}^{2}x}},n=0,1,2,\dots,x\in[0,+\infty),\label{eq:Ex2}\end{equation}
where $\beta_{n},\alpha_{n}^{1},\alpha_{n}^{2}$ are positive periodic
sequences with a period $k$. It is not necessary that $\beta_{n},\alpha_{n}^{1},\alpha_{n}^{2}$
share the same period and, in this case, we can also easily find such
a $k$. We obtain the following results by applying Theorems \ref{thm:main0}
and \ref{thm:mainthm}, respectively.
\begin{cor}
If $\Pi_{n=1}^{k}\beta_{n}\leq1$, then, for any solution of (\ref{eq:Mother_generalPielou})
with (\ref{eq:Ex2}), it holds that\[
\lim_{n\to+\infty}x_{n}=0.\]

\end{cor}
On the other hand, for the case $\Pi_{n=1}^{k}\beta_{n}>1$, we establish
the following result.
\begin{cor}
If \[
\Pi_{n=1}^{k}\beta_{n}>1,\text{ and }\Pi_{n=1}^{k}\frac{\beta_{n}}{1+\frac{\alpha_{n}^{1}}{\alpha_{n}^{2}}}<1,\]
then, there exists a $k$-periodic solution $x_{n}^{*}$ such that
$x_{n}^{*}=x_{n+k}^{*}$ which is globally attractive, that is, for
any solution of (\ref{eq:Mother_generalPielou}) with (\ref{eq:Ex2}),
it holds that \[
\lim_{n\to+\infty}(x_{n}-x_{n}^{*})=0.\]
\end{cor}
\begin{ack}
The authors are very grateful to the anonymous referee for carefully
reading and valuable comments which led to a significant improvement
of the original manuscript. This work was partially done during the
second author's stay as a member of International Research Training
Group (IGK 1529) at the Technical University Darmstadt. The authors
would like to express their gratitude to Professor Yoshiaki Muroya
for his valuable comments to this paper.
\end{ack}

\end{document}